# Consistency of Random Survival Forests

*Hemant Ishwaran and Udaya B. Kogalur*
*Cleveland Clinic*

October 27, 2018

We prove uniform consistency of Random Survival Forests (RSF), a newly introduced forest ensemble learner for analysis of right-censored survival data. Consistency is proven under general splitting rules, bootstrapping, and random selection of variables—that is, under true implementation of the methodology. A key assumption made is that all variables are factors. Although this assumes that the feature space has finite cardinality, in practice the space can be a extremely large—indeed, current computational procedures do not properly deal with this setting. An indirect consequence of this work is the introduction of new computational methodology for dealing with factors with unlimited number of labels.

## 1 Introduction

Out of the machine learning community have emerged many different types of learning algorithms. Some of these have excited tremendous interest because of their performance over benchmark data. With minimal supervision, these algorithms outperform standard methods in terms of prediction error—in some instances the difference in prediction error is so substantial there seems no way to bridge the gap. One of the most exciting algorithms to have been proposed is Random Forests (RF), an ensemble learning method introduced by Leo Breiman (Breiman, 2001). RF is an all-purpose algorithm that can be applied in a wide variety of data settings. In regression settings (i.e. where the response is continuous) the method is referred to as RF-R. In classification problems, or multiclass problems, where the response is a class label, the method is referred to as RF-C. Recently the methodology has also been extended to right-censoring survival settings, a method called random survival forests (RSF) (Ishwaran et al., 2008).

RF is considered an "ensemble learner". Ensemble learners are predictors formed by aggregating many individual learners (base learners), each of which have been constructed from different realizations of the data. There has much interest in ensemble learners because they have been shown in many instances to outperform the individual learners they are constructed from—although why this happens is not yet fully understood. A widely used ensemble technique is bagging (Breiman, 1996a). In bagging, the ensemble is formed by aggregating a base learner over independent bootstrap samples of the original data. If the base learner is a classification tree, then a classification tree is fit to each bootstrap sample, and the ensemble classifier is defined by taking a majority vote over the individual classifiers. If the base learner is a regression tree, the ensemble is the averaged tree predictor. Breiman showed that the improved performance of a bagged predictor is related to instability (Breiman, 1996a,b). If the base learner is unstable, with low bias, then the ensemble will have better performance. However, bagging a stable learner can sometimes degrade performance.

Although there are many variants of RF (Amit and Geman, 1997; Dietterich, 2000; Cutler and Zhao, 2001), the most popular, and the one we focus on here, is that described by Breiman in his software manual (Breiman, 2003). This algorithm was also discussed in Breiman (2001) under the name Forest-RI for RF random input selection, and it is also the algorithm implemented in the R-software packages `randomForest`(Liaw and Wiener, 2002), and `randomSurvivalForest` (Ishwaran and Kogalur, 2007). In this version, RF can be viewed as an extension of bagging. Using independent bootstrap samples, a random tree is grown by randomly selecting a subset of variables (features) to be used as candidate variables for splitting each node. The forest ensemble is constructed by aggregating over the random trees. The extra randomization introduced in the tree growing process is the crucial step distinguishing forests from bagging. Unlike bagging, each bootstrap tree is constructed using different variables, and not all variables are used. This is designed to encourage independence among trees, and unlike bagging, it not only reduces variance, but also bias. While this extra step might seem harmless, results using benchmark data have shown that prediction error for RF can be substantially better than bagging. In fact, performance of RF has been found comparable to other state-of-the-art methods such as boosting (Freund and Shapire, 1996), and support vector machines (Cortes and Vapnik, 1995).

Why RF works is still somewhat of a mystery. Developing theory is difficult because although each step used in RF is fairly simple, when combined, they result in an algorithm that is hard to analyze analytically. In his seminal 2001 paper (Breiman, 2001), Breiman discussed bounds on the generalization error for a forest as a trade-off involving number of variables randomly selected as candidates for splitting, and the correlation between trees. He showed as number of variables increases, strength of a tree (accuracy) increases, but at the price of increasing correlation among trees, which degrades overall performance. It is unclear, however, how tight these bounds are—Breiman himself noted they might be loose. In Lin and Jeon (2006), lower bounds for the mean-squared error for a regression forest were derived under random splitting by drawing analogies between forests and nearest neighbor classifiers. Recently, Meinshausen (2006) proved consistency of RF-R for quantile regression, and Biau, Devroye, and Lugosi (2008) proved consistency of RF-C under the assumption of random splitting.

In this paper, we prove uniform consistency of RSF. As this is a new extension of RF to right-censored survival settings not much is known about its properties. Even consistency results for survival trees are sparse in the literature. For right-censored survival data, LeBlanc and Crowley (1993) showed survival tree cumulative hazard functions are consistent for smoothed cumulative hazard functions. The method of proof used convergence results for recursive-partitioned regression trees for uncensored data Breiman et.al (1984).

We take a different approach and establish consistency by drawing upon counting process theory. We first prove uniform consistency of survival trees, and from this, by making use of bootstrap theory, we prove consistency of RSF (Section 4). These results apply to general tree splitting rules (not just random ones) and to true implementations of RSF. We make only one important assumption: that the feature space is a finite (but very large) discrete space and that all variables are factors (Section 3). In this regard we deviate from other proofs of forest consistency. These proofs all assume that the feature space is continuous; but this is problematic for two reasons. First, it requires strong assumptions about the splitting rule which may not be reflected in practice. Secondly, data in most applied statistical problems, especially those seen in medical settings, contain a mixture of both continuous and categorical data. Because node splitting for categorical variables is fundamentally different than for continuous variables, theory proven under a continuous feature space paradigm do not directly apply to most data settings.

On the other hand, to be fair to these approaches, continuous variables are often encountered in practice. Thus, it is natural to wonder if an assumption of a discrete feature space limits our results. We show in Section 5 that embedding forests in a discrete setting is realistic in that one can analyze problems with continuous variables by treating them as factors having a large number of factor labels. Indirectly, this addresses an unresolved issue related to forests and trees: namely, how to split factors with a large number of labels. We introduce new computational methodology for addressing unlimited number of labels for factors. For the interested user, we note that all computations given in the paper can be implemented using the freely available R-software package, `randomSurvivalForest` (Ishwaran and Kogalur, 2007, 2008).

## 2 Random survival forests algorithm

We begin with a high-level description of the RSF algorithm. Specific details follow.

1. Draw $B$ independent bootstrap samples from the learning data and grow a binary recursive survival tree to each bootstrap sample.

2. When growing a survival tree, at each node of the tree randomly select $p$ candidate variables to split on (use as many candidate variables as possible, up to $p - 1$, if there are less than $p$ variables available within the node). The node is split using the split that maximizes survival difference between daughter nodes (in the case of ties, a random tie breaking rule is used).

3. Grow a tree as near to saturation as possible (i.e. to full size) with the only constraint being that each terminal node should have no less than $d_0 > 0$ events.
4. Calculate the tree survival function. The forest ensemble is the averaged tree survival function.

The tree survival function calculated in Step 4 is the Kaplan-Meier (KM) estimator for the tree's terminal nodes. This can be explained formally using the following notation. Let $\mathscr{T}$ denote the terminal nodes of a survival tree, $T$. These are the extreme nodes of $T$ reached when the tree can no longer be split to form new nodes (daughters). Let $(T_{1,h}, \delta_{1,h}), \ldots, (T_{m(h),h}, \delta_{m(h),h})$ be survival times and binary $\{0, 1\}$ censoring variables for cases (individuals) in a terminal node $h \in \mathscr{T}$. An individual $i$ is said to be right-censored at time $T_{i,h}$ if $\delta_{i,h} = 0$; otherwise, if $\delta_{i,h} = 1$, the individual is said to have experienced an event at $T_{i,h}$. Let $t_{1,h} < t_{2,h} < \cdots < t_{M(h),h}$ be the $M(h)$ distinct event times. Define $d_{l,h}$ to be the number of events at time $t_{l,h}$ and $Y_{l,h}$ to be the number individuals at risk just prior to time $t_{l,h}$. The cumulative hazard function (CHF) estimate for $h$ is the Nelson-Aalen estimator

$$\hat{H}_h(t) = \sum_{t_{l,h} \leq t} \frac{d_{l,h}}{Y_{l,h}}.$$

Note that all cases within $h$ are assigned $\hat{H}_h(t)$.

For later theoretical development it will be helpful to rewrite $\hat{H}_h(t)$ using counting process notation. Let the predictable function $Y_h^{(n)}(t) = \sum_{i=1}^{m(h)} I(T_{i,h} \geq t)$ be the number of individuals in $h$ observed to be at risk just prior to $t$, and let $N_h^{(n)}(t)$ be the counting process for $h$ defined as the number of events in $[0, t]$. Define the indicator process $J_h^{(n)}(t) = I(Y_h^{(n)}(t) > 0)$. Then the Nelson-Aalen estimator for $h$ may equivalently be written as

$$\hat{H}_h(t) = \int_0^t \frac{J_h^{(n)}(s)}{Y_h^{(n)}(s)} \, dN_h^{(n)}(s),$$

where we adopt the convention that $J_h^{(n)}(s)/Y_h^{(n)}(s) = 0$ whenever $Y_h^{(n)}(s) = 0$. The KM estimator for $h$ is

$$\hat{S}_h(t) = \prod_{s \leq t} (1 - d\hat{H}_h(s)) = \prod_{s \leq t} \left(1 - \frac{dN_h^{(n)}(s)}{Y_h^{(n)}(s)}\right). \quad (1)$$

Each individual $i$ has a $d$-dimensional feature $\mathbf{x}_i$. To determine the survival function for $i$, drop $\mathbf{x}_i$ down the tree. Because of the binary nature of a survival tree, $\mathbf{x}_i$ will be assigned a unique terminal node $h \in \mathscr{T}$. The survival function for $i$ is the KM estimator for $\mathbf{x}_i$'s terminal node:

$$\hat{S}(t|\mathbf{x}_i) = \hat{S}_h(t), \quad \text{if } \mathbf{x}_i \in h.$$

Note this defines the survival function for all cases and thus defines the survival function for the tree, $T$. To make this clear, we write the tree survival function as

$$\hat{S}(t|\mathbf{x}_i) = \sum_{h \in \mathcal{T}} I(\mathbf{x}_i \in h)\hat{S}_h(t).$$

## 3 The feature space

In establishing consistency of random survival forests we assume that each coordinate $1 \leq j \leq d$ of the $d$-dimensional feature $\mathbf{X}$ is a factor (discrete nominal variable) with $1 < L_j < \infty$ distinct labels. While this assumes that the feature space $\mathcal{X}$ has finite cardinality, the actual size of $\mathcal{X}$ can be quite large, $L_1 \times \cdots \times L_d$, and moreover, the number of splits that a tree might make from such data can be even larger, depending on $d$ and $(L_j)_{j=1}^d$.

To see this, note that a split on a factor in a tree results in data points moving left and right of the parent node such that the complementary pairings define the new daughter nodes. For example, if a factor has three labels, $\{A, B, C\}$, then there are three complementary pairings (daughters) as follows: $\{A\}$ and $\{B, C\}$; $\{B\}$ and $\{C, A\}$; and $\{C\}$ and $\{A, B\}$. In general, for a factor with $L_j$ distinct labels, there are $2^{L_j-1} - 1$ distinct complementary pairs. Thus, the total number of splits evaluated when splitting the root node for a survival tree when all variables are factors can be as much as

$$\text{maximum number root-node splits} = \sum_{j=1}^{d} 2^{L_j-1} - d.$$

Following the root-node split, are splits on the resulting daugther nodes, and their daughter nodes, recursively, with each subsequent generation requiring a large number of evaluations. Each evaluation can result in a new tree, thus showing that number of trees (space of trees) associated with $\mathcal{X}$ can be extremely large.

## 4 Properties of survival forests

A reasonable criteria for consistency of a random survival forest is that the ensemble survival function converges to the population survival function. We first consider conditions needed for consistency of a survival tree. Consistency of forests are then deduced by utilizing bootstrap theory.

### 4.1 Assumptions

Let $(\mathbf{X}, T, \delta), (\mathbf{X}_1, T_1, \delta_1), \ldots, (\mathbf{X}_n, T_n, \delta_n)$ be i.i.d. random elements such that $\mathbf{X}$, the feature, takes values in $\mathcal{X}$, a discrete space as described in Section 3. Here $T =$ $\min(T^0, C)$ is the observed survival time and $\delta = I(T^o \le C)$ is the binary $\{0,1\}$ censoring value, where it is assumed that $T^o$, the true event time, is independent of $C$, the censoring time. Furthermore, it is assumed that $\mathbf{X}$ is independent of $\delta$. The collection of values $\{(\mathbf{X}_i, T_i, \delta_i)\}_{i=1}^n$ are referred to as the *learning data* and are used in the construction of the forest (recall the algorithm of Section 2). It is assumed that $(\mathbf{X}, T, \delta)$ has joint distribution $\mathbb{P}$. The marginal distribution for $\mathbf{X}$ is denoted by $\mu$ and defined via $\mu(A) = \mathbb{P}\{\mathbf{X} \in A\}$ for all subsets $A$ of $\mathscr{X}$. It is assumed that $\mu(A) > 0$ for each $A \ne \emptyset$.

The true survival function, or population parameter, is assumed to be of the form

$$S(t|\mathbf{X}) := \mathbb{P}\{T^o > t|\mathbf{X}\} = \sum_{\mathbf{x} \in \mathscr{X}} I(\mathbf{X} = \mathbf{x}) \exp\left(-\int_0^t \alpha(s|\mathbf{x})\,ds\right), \quad (2)$$

where $\alpha(\cdot|\mathbf{x})$ is the non-negative hazard function for the subpopulation $\mathbf{X} = \mathbf{x}$.

### 4.2 Uniform consistency of survival trees

The following result, showing uniform consistency of a survival tree, is a consequence of the uniform consistency of the KM estimator. Let $\tau = \min\{\tau(\mathbf{x}) : \mathbf{x} \in \mathscr{X}\}$, where $\tau(\mathbf{x}) = \sup\{t : \int_0^t \alpha(s|\mathbf{x})ds < \infty\}$.

**Theorem 1.** *Let $t \in (0, \tau)$. If $\mathbb{P}\{C > 0\} > 0$, and $\alpha(\cdot|\mathbf{x})$ is strictly positive over $[0,t]$ for at least one $\mathbf{x} \in \mathscr{X}$, then*

$$\sup_{s \in [0,t]} |\mu(\hat{S}(s|\mathbf{X})) - \mu(S(s|\mathbf{X}))| \xrightarrow{\mathrm{P}} 0, \quad as\ n \to \infty.$$

*(Note that we are using the linear functional notation for expectation in the above expression).*

### 4.3 Uniform consistency of survival forests

A random survival forest is a collection of random survival trees, each grown from independent bootstrap samples of the learning data, $\mathscr{L} = \{(\mathbf{X}_i, T_i, \delta_i)\}_{i=1}^n$. Thus, in order to prove consistency of RSF, we must extend our previous results to include bootstrap resampling.

Let $\mathscr{L}^* = \{(\mathbf{X}_i^*, T_i^*, \delta_i^*)\}_{i=1}^n$ denote a bootstrap sample of the learning data. Let $T^*$ be the survival tree grown from $\mathscr{L}^*$ and let $\hat{S}^*(t|\mathbf{x})$ be the KM estimator for $T^*$:

$$\hat{S}^*(t|\mathbf{x}) = \sum_{h \in \mathscr{T}^*} I(\mathbf{x} \in h)\hat{S}_h^*(t),$$

where $\hat{S}_h^*(t)$ is defined similar to (1), and the above sum is over $\mathscr{T}^*$, the set of terminal nodes of $T^*$.

A random survival forest comprised of $B$ survival trees has an ensemble survival function

$$\hat{S}_e(t|\mathbf{x}) = \frac{1}{B} \sum_{b=1}^{B} \hat{S}_b^*(t|\mathbf{x}),$$

where $\hat{S}_b^*(t|\mathbf{x})$ is the survival function for the survival tree grown using the $b$-th bootstrap sample. We show consistency of RSF by establishing consistency of $\mu(\hat{S}_b^*(t|\mathbf{X}))$ for each $b$ (since then consistency of the ensemble, $\mu(\hat{S}_e(t|\mathbf{X}))$, holds automatically).

**Theorem 2.** *Let $\tau^* = \min(\tau, \sup(F))$, where $\sup(F)$ is the upper limit of the support of $F(s) = 1 - \mathbb{P}\{T^o > s\}\mathbb{P}\{C > s\}$. Then under the same conditions as in Theorem 1, for each $t \in (0, \tau^*)$:*

$$\sup_{s \in [0,t]} |\mu(\hat{S}^*(s|\mathbf{X})) - \mu(S(s|\mathbf{X}))| = o_p^*(1) + o_p(1),$$

*where $o_p^*$ stands for $o_p$ in bootstrap probability for almost all $\mathscr{L}$-sample sequences; i.e. with probability one under $\mathbb{P}^\infty$.*

## 4.4 Uniform approximation by forests

Theorem 2 establishes consistency of a bootstrapped survival tree, and from this consistency of a survival forest follows. While this is a useful line of attack for establishing large sample properties of forests, it does not convey how in practice a forest might improve inference over a single tree. Indeed, in finite sample settings, a forest of trees can have a decided advantage when approximating the true survival function.

To show this we make use of the following, somewhat idealized, setting. Suppose that for each $b$ we are allowed to construct a binary survival tree $T_b$ from a prechosen learning data set $\mathscr{L}_b = \{(\mathbf{X}_{b,i}, T_{b,i}, \delta_{b,i})\}_{i=1}^n$ in any manner we choose. The only constraint being that each terminal node of $T_b$ must contain at least $d_0 = 1$ events. Let $S_b(t|\mathbf{x})$ be the KM tree survival function for $T_b$, and let

$$S_e(t|\mathbf{x}) = \sum_{b=1}^{B} W_b S_b(t|\mathbf{x}) \qquad (3)$$

be the ensemble survival function for the forest comprising $(T_b)_{b=1}^B$, where $W_b \geq 0$ are forest weights. The next theorem shows that one can always find an ensemble that uniformly approximates the true survival function (2). Trees do not possess this property.

**Theorem 3.** *If $n > d$, and $s \in [0, \tau)$, then for each $\varepsilon > 0$ there exists an ensemble survival function (3) for a survival forest comprised of $B = B(\varepsilon)$ survival trees, with each tree consisting of $d + 1$ terminal nodes, such that*

$$\int_0^s \int_{\mathscr{X}} (S_e(t|\mathbf{x}) - S(t|\mathbf{x}))^2 \, \mu(d\mathbf{x}) \, dt \leq \varepsilon.$$

## 5 Empirical results

Our theory has been based on the assumption that all $x$-variables are factors, but in practice one often encounters data with continuous variables. Here we show that one can discretize continuous variables and treat them as factors without unduly affecting prediction error and inference: thus showing our theory can be extrapolated to general data settings. At the same time, an indirect consequence of this work is the introduction of new computational methodology for efficient splitting of factors and for dealing with factors with unlimited number of labels.

For illustration we consider the primary biliary cirrhosis (PBC) data of Fleming and Harrington (1991). The data is from a randomized clinical trial studing the effectiveness of the drug D-penicillamine on PBC. The dataset involves 312 individuals and contains 17 variables as well as censoring information and time until death for each individual. Of the 17 features, seven are discrete and 10 are continuous. Each of the 10 continuous variables were discretized and converted to a factor with $L$ labels. We investigated different amounts of granularity: $L = 2, \ldots, 30$.

For each level of granularity, $L$, we fit a survival forest of 1000 survival trees using log-rank splitting with node-adaptive random splits. Splits for nodes were implemented as follows. A maximum of "nsplit" complementary pairs were chosen randomly for each of the $p$ randomly selected candidate variables within a node (if nsplit exceeded the number of cases in a node, then nsplit was set to the size of the node). Log-rank splitting was applied to the randomly selected complementary pairs and the node was split on that variable and complementary pair maximizing the log-rank test. Five different values for nsplit were tried: nsplit= $5, 10, 20, 50, 1024$. All computations were implemented using the R-software package, `randomSurvivalForest` (Ishwaran and Kogalur, 2007, 2008).

The top plot in Figure 1 shows out-of-bag prediction error as a function of granularity and nsplit value. As granularity rises, prediction error increases—but this increase is reasonably slow and well contained with larger values of nsplit. This is quite remarkable because total number of complementary pairs with a granularity level of $L = 30$ is on order $2^{30}$ (over 1 billion pairs) and yet our results show that using only 50 randomly selected complementary pairs keeps prediction error in check.

Prediction error measures overall performance, but we should also consider how inference for a variable is affected by increasing granularity. To study this we looked at variable importance (VIMP). VIMP measures predictiveness of a variable, adjusting for all other variables. Positive values of VIMP indicating predictiveness, and negative and zero values indicating noise (Ishwaran et al., 2008). For each forest used in Figure 1 we dropped bootstrapped data down the forest and computed VIMP for each variable. This was repeated 1000 times independently resulting in a bootstrap distribution for VIMP. The bottom plot of Figure 1 displays the 68% bootstrap confidence region from this distribution. The analysis was restricted to only those forests grown under an nsplit value of 1024 but was carried out for each level of granularity as in Figure 1 (color coding scheme used to depict granularity is described in the caption of the figure). Overall, one can see that the bootstrap confidence regions are relatively robust to the level of granularity.

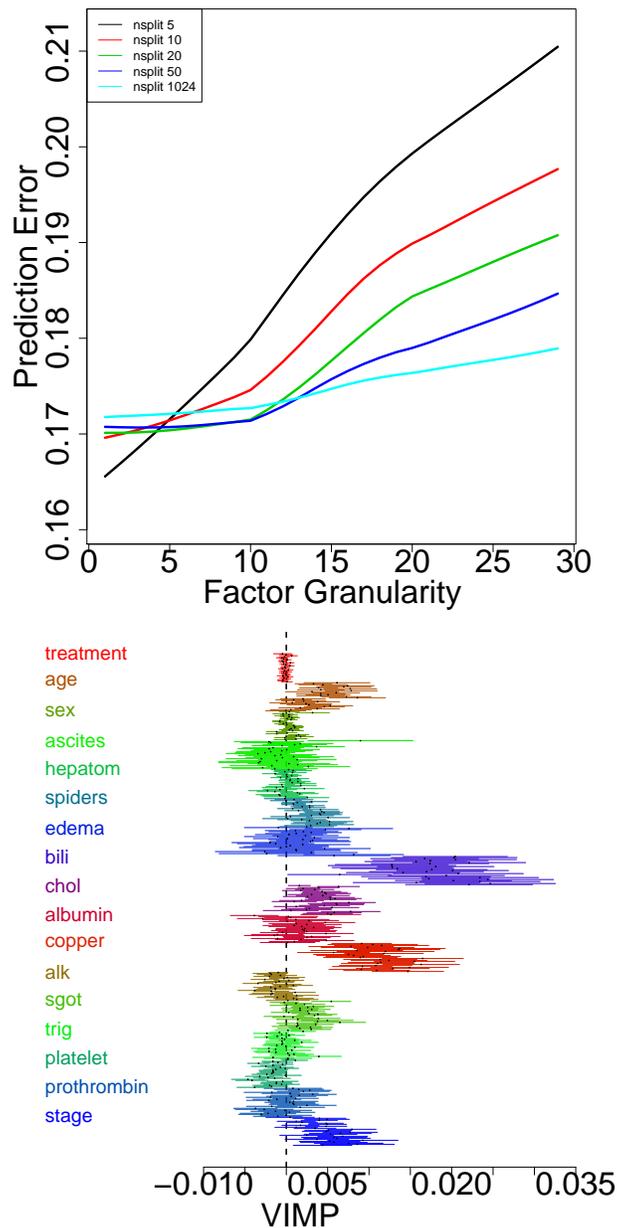

Figure 1: RSF analysis of PBC data using 1000 trees with random log-rank splitting where variables, both nominal and continuous, were discretized to have a maximum number of labels (factor granularity). Top figure is out-of-bag prediction error versus factor granularity, stratified by number of random splits used for a node, nsplit. Bottom figure shows 68% bootstrap confidence region for variable importance (VIMP) from 1000 bootstrap samples using an nsplit value of 1024 for each factor granularity value in the top figure. Color coding is such that the same color has been used for a variable over the different granularity values (factor granularity for a variable increases going from top to bottom).

### 5.1 Remarks

1. Not only does random splitting maintain good prediction error performance, but it may actually help mitigate node-splitting bias. It has been noted in the literature that forests have a tendency to favor continuous variables and factors with a large number of labels because of a type-I error effect—that is, with more values to split on, there is an increased probability of finding a spurious effect. To investigate the effect to which random splitting can alleviate such bias we expanded the PBC data to include 50 independent noise variables. 25 of these were randomly simulated from a standard normal distribution, thus representing continuous variables, the other 25 were discrete variables, randomly simulated from a two-point distribution having equal probability for each class.

   The data was discretized and the bootstrap VIMP distribution for each variable calculated as in the previous section (in total there was 67 variables). As in the previous section granularity levels of $L = 2, \ldots, 30$ were investigated. The analysis was restricted to an nsplit value of 1024.

   The 68% bootstrap VIMP confidence regions are depicted in Figure 2 (results are displayed similar to the bottom plot of Figure 1). Clearly VIMP distributions for continuous noise variables are wider than that for discrete noise variables. However, nearly all distributions for continuous noise variables contain zero, even for high levels of granularity, clearly showing that random splitting is helping to mitigate selection bias.

2. Variable selection by purposefully adding noise variables has been succesfully applied to certain regression models (Wu, Boos and Stefasnki, 2007) and the results from Figure 2 suggest a similar idea may work for forests. One implemention would be to introduce a reasonably large number of noise variables and use the combined bootstrap distribution for the variables to determine a VIMP threshold (combining the distributions should yield a more stable threshold value). A nice feature of forests is that the inclusion of even a fairly large number of noise variables should not impact prediction error or the VIMP of non-noise variables.

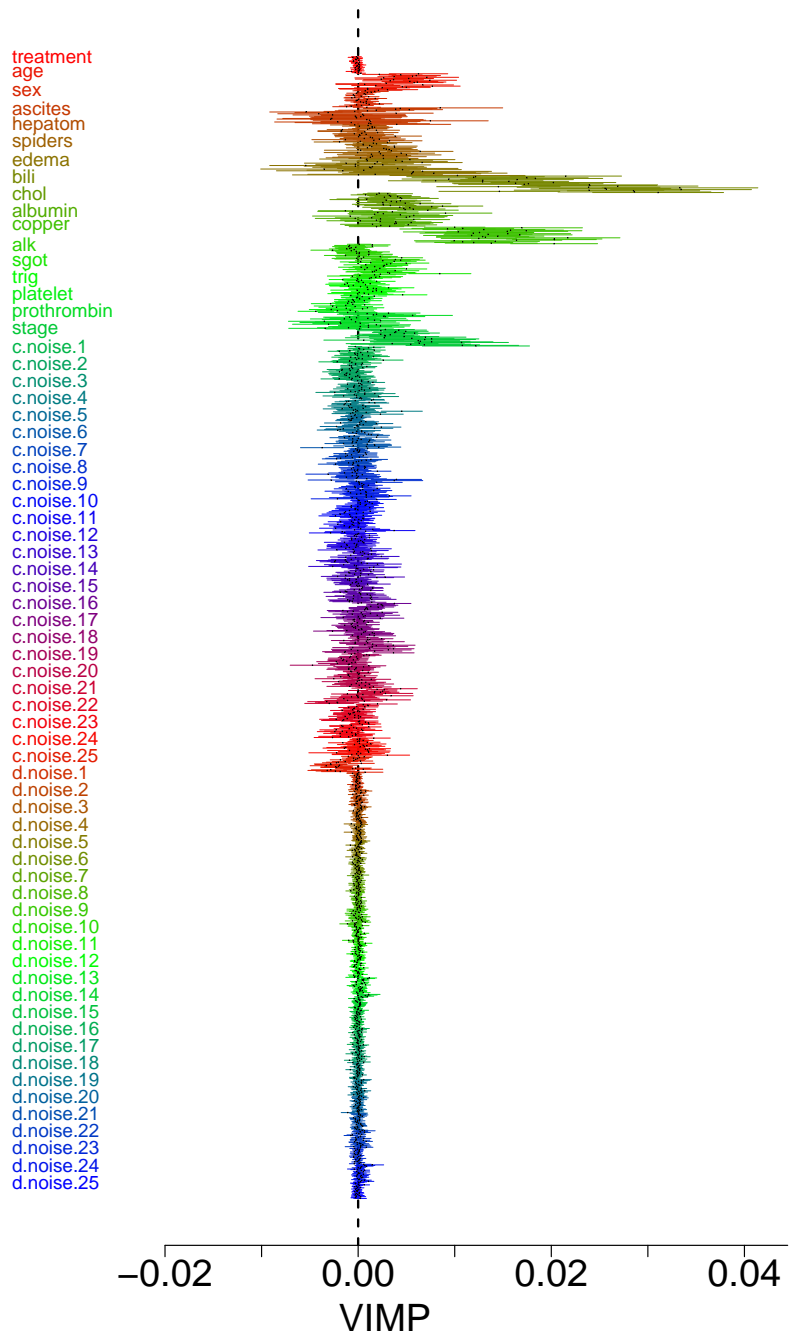

Figure 2: Bootstrapped VIMP from PBC analysis where data has been expanded to include 25 continuous and 25 discrete binary noise variables. Discrete noise variables encoded using a "d", whereas continuous noise variables have codes starting with "c". VIMP calculated and displayed as in bottom plot of Figure 1.

3. Although our examples considered factors with no more than 30 labels, we are able to implement splitting on factors with unlimited number of labels in `randomSurvivalForest`. Here is a brief description of this methodology.

   The most basic issue is how to represent a split on a factor. Here we emulate `randomForest`, the R-software package based on Breiman and Cutler and ported to R by Liaw and Wiener (Liaw and Wiener, 2002). The strategy is to immutably map the labels in the factor to the bit positions of a 32-bit integer. A split is then uniquely defined by moving all labels corresponding to bits ON (equal to one) to the left daughter and moving the rest, corresponding to bits OFF (equal to zero) to the right daughter. Note that the left daughter and right daughter define complementary subsets of the factor labels.

   For deterministic splits on a factor with 32 labels, all possible complementary pairs of subsets are enumerated. In general there are $2^L - 1$ such pairs. It is clear that, with $L = 32$ labels, enumerating the number of complementary pairs explicitly becomes memory intensive. In fact, on an architecture in which an unsigned integer is a 32-bit word, the maximum value representable is a factor with 32 labels.

   Deterministic splitting requires that we construct the binary representation of all the possible complementary pairs in a factor—but this is memory intensive. Our solution is to allow factors with unlimited labels, but to restrict factors with labels greater than 32 labels to node-adaptive random splitting (with nsplit set to the cardinality of the working node). We thus completely avoid the overhead of enumerating the complementary pairs by constructing and discarding each complementary pair individually.

   Once a complementary pair is identified for splitting a tree node, representing the split point in its binary complementary pair format is an extension of the single word case. We simply use an array of 32-bit unsigned integers that is sufficient in length to represent all the labels in the factor. For example, a factor with 512 labels requires a vector composed of 16 unsigned integers. We refer to this vector as a multi-word complementary pair (MWCP). By a simple book-keeping mechanism, we are able to store the vector and its length, and are thus able to recover the split point for later use (for example for prediction on test data).

## 6 Proofs

*Proof of Theorem 1.* Independence of $T^o$ and $C$ ensures that

$$\mathbb{P}\{\delta = 1\} = \mathbb{P}\Big(\mathbb{P}\{T^o \leq c | C = c\}\Big) \geq \mathbb{P}\{T^o \leq c\}\mathbb{P}\{C > c\}, \quad \text{for any } c > 0.$$

The assumption $\mathbb{P}\{C > 0\} > 0$ implies that the censoring distribution has mass bounded away from the origin. Thus, we can find a $c > 0$ such that $\mathbb{P}\{C > c\} > 0$. The assumed form of the survival function (2) ensures that the distribution function for $T^0$ is continuous over $[0, t]$. Combining this with the assumption $\alpha(\cdot|\mathbf{x})$ is strictly positive for some $\mathbf{x}$, implies that $\mathbb{P}\{T^o \leq c\} > 0$, and hence, $\mathbb{P}\{\delta = 1\} > 0$. Recall that a survival tree is grown to full length with the proviso that a terminal node should have no less than $d_0 > 0$ events. Let $A \subseteq \mathscr{X}$ be any non-null set. Then by the law of large numbers, and by the assumed independence of $\mathbf{X}$ and $\delta$,

$$\frac{1}{n}\sum_{i=1}^n I(\mathbf{X}_i \in A, \delta_i = 1) \overset{\text{a.s.}}{\to} \mathbb{P}\{\mathbf{X} \in A, \delta = 1\} = \mu(A)\mathbb{P}\{\delta = 1\} > 0.$$

Therefore,

$$I\left(\sum_{i=1}^n I(\mathbf{X}_i \in A, \delta_i = 1) \geq d_0\right) \overset{\text{a.s.}}{\to} 1. \tag{4}$$

Thus, the constraint that a terminal node must have $d_0$ events holds almost surely if the terminal node is $A$. Furthermore, because the survival tree is grown to full length, this shows that $A$ must be a distinct value $\mathbf{x} \in \mathscr{X}$. If it were not, then this would imply the tree stopped splitting at a node comprised of more than one $\mathbf{x}$, because the requirement of $d_0 > 0$ events could not be met under any split. That is, any split on this node yields daughters $A_1$ and $A_2$ with at least one daughter having fewer than $d_0$ deaths. But this contradicts (4) which holds for any $A$. Thus, the tree almost surely splits on all possible values of $\mathscr{X}$ and has terminal nodes for each distinct $\mathbf{x} \in \mathscr{X}$. In other words,

$$\hat{S}(s|\mathbf{X}) = \sum_{\mathbf{x} \in \mathscr{X}} I(\mathbf{X} = \mathbf{x} = h)\hat{S}_h(s) + o_p(1).$$

Note that the $o_p(1)$ term is uniform in $s$.

Let $Y^{(n)}(s|\mathbf{x}) = \sum_{i=1}^n I(T_i \geq s, \mathbf{X}_i = \mathbf{x})$ be the number of cases with feature $\mathbf{x}$ who are at risk just prior to $s$. Define $J^{(n)}(s|\mathbf{x}) = I(Y^{(n)}(s|\mathbf{x}) > 0)$. Now if we can show that for each $\mathbf{x} \in \mathscr{X}$, as $n \to \infty$,

$$\int_0^t \frac{J^{(n)}(s|\mathbf{x})}{Y^{(n)}(s|\mathbf{x})}\alpha(s|\mathbf{x})\, ds \overset{\text{P}}{\to} 0, \tag{5}$$

and

$$\int_0^t (1 - J^{(n)}(s|\mathbf{x}))\alpha(s|\mathbf{x})\, ds \overset{\text{P}}{\to} 0, \tag{6}$$

then by Theorem IV.3.1 (Andersen et al., 1993), for each $h = \mathbf{x}$:

$$\sup_{s\in[0,t]} |\hat{S}_h(s) - S(s|\mathbf{x})| \xrightarrow{\text{P}} 0, \quad \text{as } n \to \infty.$$

This would establish the result, because

$$\sup_{s\in[0,t]} |\mu(\hat{S}(s|\mathbf{X})) - \mu(S(s|\mathbf{X})|$$
$$\leq \sum_{h=\mathbf{x}} \mu\{\mathbf{X}=\mathbf{x}\} \sup_{s\in[0,t]} |\hat{S}_h(s) - S(s|\mathbf{x})| + o_p(1).$$

Therefore, to complete the proof we need only to verify that conditions (5) and (6) hold. By the definition of $\tau$, $\sup_{s\in[0,t]} \alpha(s|\mathbf{x}) < \infty$ and thus a sufficient condition for (5) and (6) is that

$$\inf_{s\in[0,t]} Y^{(n)}(s|\mathbf{x}) \xrightarrow{\text{P}} \infty.$$

This condition holds by noting that for each $s \in [0,t]$,

$$n^{-1}Y^{(n)}(s|\mathbf{x}) \geq \frac{1}{n}\sum_{i=1}^{n} I(T_i \geq t, \delta_i = 1, \mathbf{X}_i = \mathbf{x})$$
$$\xrightarrow{\text{a.s.}} \mu(\mathbf{X}=\mathbf{x})\mathbb{P}\{T^o \geq t|\mathbf{x}\} > 0.$$

Note that $\mathbb{P}\{T^o \geq t|\mathbf{x}\} > 0$ because of (2) and the definition of $\tau$. □

*Proof of Theorem 2.* Let $(M^*_{n,1}, \ldots, M^*_{n,n})^T$ be a multinomial random vector from $n$ trials in which each cell has probability $1/n$ (as customary we use a "*" to indicate that randomness comes from bootstrapping). For each non-null $A \subseteq \mathscr{X}$,

$$\frac{1}{n}\sum_{i=1}^{n} I(\mathbf{X}_i^* \in A, \delta_i^* = 1)$$
$$\stackrel{d^*}{=} \frac{1}{n}\sum_{i=1}^{n} I(\mathbf{X}_i \in A, \delta_i = 1)M^*_{n,i}$$
$$= \frac{1}{n}\sum_{i=1}^{n} I(\mathbf{X}_i \in A, \delta_i = 1) + \frac{1}{n}\sum_{i=1}^{n} I(\mathbf{X}_i \in A, \delta_i = 1)(M^*_{n,i} - 1)$$
$$= \mathbb{P}\{\mathbf{X} \in A, \delta = 1\} + o_p^*(1), \quad \text{a.s.,}$$

where the almost sure convergence is in $\mathbb{P}$-probability and the $o_p^*(1)$ term follows

from

$$\mathbb{P}^*\left\{\frac{1}{n}\sum_{i=1}^n I(\mathbf{X}_i \in A, \delta_i = 1)(M_{n,i}^* - 1) \geq \varepsilon\right\}$$

$$\leq \frac{1}{\varepsilon^2 n^2}\mathbb{E}^*\left(\sum_{i=1}^n I(\mathbf{X}_i \in A, \delta_i = 1)(M_{n,i}^* - 1)\right)^2$$

$$\leq \frac{1}{\varepsilon^2 n^2}\sum_{i=1}^n \mathrm{Var}^*(M_{n,i}^*) + \frac{1}{\varepsilon^2 n^2}\sum_{i\neq j}|\mathrm{Cov}^*(M_{n,i}^*, M_{n,j}^*)|$$

$$= O(n^{-1}).$$

In the proof of Theorem 1 it was shown that $\mathbb{P}\{\mathbf{X} \in A, \delta = 1\} > 0$. From this, and ussing similar arguments as in that proof, it follows that

$$\hat{S}^*(s|\mathbf{X}) = \sum_{\mathbf{x}\in\mathscr{X}} I(\mathbf{X} = \mathbf{x} = h)\hat{S}_h^*(s) + o_p^*(1),$$

where the $o_p^*(1)$ term is uniform in $s$. Notice that

$$\mu(\hat{S}^*(s|\mathbf{X})) - \mu(S(s|\mathbf{X}))$$
$$= \left[\mu(\hat{S}^*(s|\mathbf{X})) - \mu(\hat{S}(s|\mathbf{X}))\right] + \left[\mu(\hat{S}(s|\mathbf{X})) - \mu(S(s|\mathbf{X}))\right].$$

The second term in square brackets is $o_p(1)$ uniformly in $s$ by Theorem 1. To deal with the first term, we use the representation for $\hat{S}(s|\mathbf{X})$ given in the proof of Theorem 1, to obtain

$$\mu(\hat{S}^*(s|\mathbf{X})) - \mu(\hat{S}(s|\mathbf{X}))$$
$$= \sum_{h=\mathbf{x}}\mu(\mathbf{X} = \mathbf{x})\left(\hat{S}_h^*(s) - \hat{S}_h(s)\right) + o_p^*(1) + o_p(1).$$

A bootstrap sample of $\mathscr{L}$ can be drawn equivalently using a two-stage process by first drawing a multinomial vector $(n_h^*)_h$ from $n$ trials with each cell $h = \mathbf{x}$ having probability $\mu(\mathbf{X} = \mathbf{x} = h)$ and then drawing a bootstrap sample of size $n_h^*$, for each $h$, from $\mathscr{L}_h = \{(\mathbf{X}_i, T_i, \delta_i) : \mathbf{X}_i = \mathbf{x} = h\}$. It is not hard to show that $n_h^*/n_h \xrightarrow{p^*} 1$, where $n_h = |\mathscr{L}_h| = \sum_{i=1}^n I(\mathbf{X}_i = \mathbf{x} = h) \xrightarrow{a.s.} \infty$. Thus, to complete the proof it suffices to show convergence of $\hat{S}_h^*(s) - \hat{S}_h(s)$, for each $h = \mathbf{x}$, where $\hat{S}_h^*(s)$ is derived from a bootstrap sample of size $n_h$ from $\mathscr{L}_h$. We use part (b) of Lemma 3 from Lo and Singh (1985) for this. This theorem applies if $s < \sup(F)$ under the random censorship model for a continuous survival distribution. All these conditions are met under our assumptions. Thus applying this Lemma, one can show that

$$\sup_{s\in[0,t]}|\hat{S}_h^*(s) - \hat{S}_h(s)| = O_p^*(n_h^{-1/2}\log n_h^{1/2}) = o_p^*(1).$$

□

*Proof of Theorem 3.* It suffices to show that for a given $\mathbf{x}$, we can find an ensemble $S_e(\cdot|\mathbf{x}')$ that is zero if $\mathbf{x}' \neq \mathbf{x}$, and for $\mathbf{x}' = \mathbf{x}$, uniformly approximates $S(\cdot|\mathbf{x})$ over $[0, s]$ (this suffices because we can always combine such ensembles to uniformly approximate $S_e(\cdot|\mathbf{x})$ for all $\mathbf{x}$). Choose $\mathscr{L}_b$ such that $\delta_{b,i} = 1$ for all $b$ and $i$. By repeated splitting on the left, it is clear we can use $d$ splits to construct a tree $T_b$ having $d+1$ terminal nodes such that the left-most daughter node corresponds to $\mathbf{x}$. Because $n > d$, we can assign at least one event to the left-most node. For concreteness, assume that the node contains exactly one event, with event time $T > 0$. Over the remaining $d$ terminal nodes assign event times $T = 0$ for all cases. Thus $S_b(\cdot|\mathbf{x}') = 0$ if $\mathbf{x}' \neq \mathbf{x}$. On the other hand, if $\mathbf{x}' = \mathbf{x}$, then $S_b(t|\mathbf{x})$ is a step function with value 1 if $t < T$ and value 0 if $t \geq T$. Because $S(\cdot|\mathbf{x})$ is continuous [by definition (2)], it is uniformly continuous over the compact set $[0, s]$. A uniformly continuous, monotonically decreasing function over a compact set can be uniformly approximated by a linear combination of a finite number of step functions such as $S_b(\cdot|\mathbf{x})$. Thus one can construct a finite number of survival trees like $T_b$, that when suitably weighted, uniformly approximates $S(\cdot|\mathbf{x})$ over $[0, s]$. □